\documentclass[10pt]{amsart}
\usepackage{amssymb,amsmath,amscd}
\usepackage[T1]{fontenc}

\usepackage{color}
\usepackage{epsfig}

\usepackage[all]{xy}

\newtheorem{theorem}{Theorem}[section]

\newtheorem{proposition}[theorem]{Proposition}

\newcommand{\Z}{\mathbb{Z}}

\newcommand{\factor}[2]{\raisebox{0.1ex}{$#1$}\raisebox{-0.3ex}{/}\raisebox{-0.9ex}{$#2$}}

\newcounter{Lcount}

\begin{document}

\title[Computing the Cohomology of $Sub_3 S^n$]{On the calculation of the cohomology of the third finite subset space of spheres}
\author{Simon C. F. Rose}
\address{Department of Mathematics, University of British Columbia, Canada}
\email{scfr@math.ubc.ca}



\bibliographystyle{amsplain}

\maketitle

\begin{abstract}
In this paper we provide a computation of the mod 2 cohomology groups of the third finite subset space of the sphere $S^n$ using known results about the cohomology of the symmetric product of spheres.
\end{abstract}

\section{Introduction}

The space of all non-empty finite subsets (of cardinality at most $k$) of a given space $X$, denoted at times $\exp_k(X), Sub_k(X), \mathcal{F}_k(X)$ or $X^{(k)}$, has been studied since 1949 \cite{borsuk}. The most celebrated result is in a correction to this initial paper wherein it is proven \cite{bott} that $Sub_3S^1 \cong S^3$.  In addition, a result due to E. Shchepin (proven in \cite{mostovoy, msc_thesis}) shows that the inclusion $S^1 \cong Sub_1 S^1 \hookrightarrow Sub_3 S^1 \cong S^3$ is in fact the trefoil knot.

More recently (see \cite{tuffley-2002}), the result that $Sub_3 S^1 \cong S^3$ has been extended to show that $Sub_k(S^1)$ has the homotopy type of an odd dimensional sphere for all $k$.

In this paper we extend existing knowledge in the other direction---that is, we begin studying the spaces $Sub_3 S^n$ for higher values of $n$.

For the duration of this paper all coefficients will be assumed to be $\factor{\Z}{2}$, and the following notation will be used. $C_k(X)$ will denote the unordered configuration space of points, while $Sym^k(X)$ will denote the $k$-th symmetric power of $X$. We note of course that $Sym^2(X)$ is homeomorphic to $Sub_2(X)$ for all $X$.

The main theorem of this paper is the following.

\begin{theorem}\label{thm_main}
The cohomology ring $H^*(Sub_3S^n)$ is given by
\[
H^k(Sub_3S^n) =
\begin{cases}
\factor{\Z}{2} & k = 0, 2n+1, 2n+2, \ldots, 3n \\
\factor{\Z}{2} \oplus \factor{\Z}{2} & k = n + 2, \ldots, 2n \\
0 & otherwise
\end{cases}
\]
with trivial product structure.
\end{theorem}

\section{Proof of Theorem \ref{thm_main}}

We need the following preliminary results, both from the paper \cite{nakaoka}.

\begin{proposition}\label{prop_coho_exp_2}
The cohomology groups $H^k(Sub_2S^n)$ are given by
\[
H^k(Sub_2S^n) =
\begin{cases}
\factor{\Z}{2} & k = n, n + 2, n + 3, \ldots, 2n -1 , 2n \\
0 & otherwise
\end{cases}.
\]
Moreover, if $e_k$ is the generator of $H^k(Sub_2S^n)$ (for $k = n, n + 2, n + 3, \ldots, 2n$), then $Sq^ie_n = e_{n+i}$ for $i = 2, \ldots, n$.
\end{proposition}

\begin{proposition}\label{prop_coho_sym_3}
The cohomology groups $H^k(SP^3 S^n)$ are given by
\[
H^k(SP^3 S^n) =
\begin{cases}
\factor{\Z}{2} & k = n, n + 2, n + 3, \ldots, 2n -1 , 2n, 2n + 2, 2n + 3, \ldots, 3n \\
0 & otherwise
\end{cases}.
\]
Moreover, if $e_k$ is the generator of $H^k(SP^3 S^n)$, then we have the following relations.
\[
e^n \smile e^n = e^{2n} \qquad e^n \smile e^n \smile e^n = e^{3n}
\]
and, for $2 \leq i \leq n$,
\[
Sq^ie^n = e_{n+i}.
\]
All other products and squares are zero.
\end{proposition}

We now prove our main theorem.

\begin{proof}[Proof of Theorem \ref{thm_main}]
We begin with the following pushout diagram.
\begin{equation}\label{eq_pushout}
\xymatrix{
S^n \times S^n \ar[r]^{\alpha}\ar[d]_{q} & SP^3S^n \ar[d]^f \\
Sub_2 S^n \ar@{^{(}->}[r] & Sub_3 S^n
}
\end{equation}
where the map $q$ is simply the quotient map, and the map $\alpha$ takes $(x, y)$ to $[x,x,y]$. It is easy then to see that this is a pushout as claimed. From this we obtain a Mayer-Vietoris type sequence
\begin{equation}\label{eq_mv}
\cdots \to H^k(Sub_3 S^n) \to H^k(Sub_2 S^n) \oplus H^k(SP^3 S^n) \to H^k(S^n \times S^n) \to \cdots.
\end{equation}

Since for most values of $k$ we have that $H^k(S^n \times S^n)$ is zero it follows that
\[
H^k(Sub_3 S^n) \cong H^k(Sub_2 S^n) \oplus H^k(SP^3 S^n)
\]
for $k \neq n, n + 1, 2n, 2n + 1$ where we have to be somewhat more careful.

The two portions of the exact sequence are then
\[
0 \to H^n(Sub_3 S^n) \to \factor{\Z}{2} \oplus \factor{\Z}{2} \to \factor{\Z}{2} \oplus \factor{\Z}{2} \to H^{n+1}(Sub_3 S^n) \to 0
\]
and
\begin{multline*}
0 \to H^{2n}(Sub_3 S^n) \to H^{2n}(\exp_2S^n)\oplus H^{2n}(SP^3S^n) \to \\
\to H^{2n}(S^n \times S^n) \to H^{2n+1}(Sub_3 S^n) \to H^{2n+1}(SP^3S^n) \to 0
\end{multline*}

We note from \cite{KS} that $Sub_3 X$ is $(r + 1)$-connected when $X$ is $r$-connected. Here this yields that $Sub_3 S^n$ is $n$-connected, and so in particular the first sequence shows that we have $H^n(Sub_3 S^n)$ and $H^{n+1}(Sub_3 S^n)$ both zero.

We next look at the commuting squares
\[
\xymatrix{
H^n(Sub_2 S^n) \ar[r]\ar[d]^{Sq^n}_{\cong} & H^n(S^n \times S^n) \ar[d]^{Sq^n} \\
H^{2n}(Sub_2 S^n) \ar[r] & H^{2n}(S^n \times S^n)
}
\]
and
\[
\xymatrix{
H^n(SP^3S^n) \ar[r]\ar[d]^{Sq^n}_{\cong} & H^n(S^n \times S^n) \ar[d]^{Sq^n} \\
H^{2n}(SP^3S^n) \ar[r] & H^{2n}(S^n \times S^n)
}
\]
Since the squares act trivially on $H^*(S^n \times S^n)$, the composition around the right must be zero. This can only be true if the maps in dimension $2n$ are also be zero. The conclusion then follows.

To determine the product structure, we note that the only non-trivial products that could occur are from elements $\alpha_i \in H^{n+k_i}(Sub_3 S^n)$. However, we note that if $\alpha_1 \smile \alpha_2 \neq 0$ then $f^*(\alpha_1) \smile f^*(\alpha_2) = f^*(\alpha_1 \smile \alpha_2) \neq 0$. However, from Proposition \ref{prop_coho_sym_3} we see that this cannot be true.
\end{proof}

\section{Further work}

Ideally we would like to pursue this argument further to compute the cohomology of $Sub_k S^n$ for higher values of $k$ using a similar argument. We note that we always have an analagous pushout to \eqref{eq_pushout}, namely
\begin{equation}\label{eq_pushout_k}
\xymatrix{
(S^n)^{\times k} \ar[r]^{\alpha}\ar[d]^{q} & SP^{k+1}S^n \ar[d] \\
\exp_kS^n \ar@{^{(}->}[r] & \exp_{k+1}S^n
}
\end{equation}
where, again $q$ is the quotient map and $\alpha(x_1, \ldots, x_k) = [x_1, x_1, x_2, \ldots, x_k]$. However, the Mayer-Vietoris type sequence of \eqref{eq_mv} no longer holds, and so we must replace \eqref{eq_pushout_k} with
\begin{equation}\label{eq_pushout_k2}
\xymatrix{
\Delta^{k+1} \ar@{^{(}->}[r]\ar[d]^{q} & SP^{k+1}S^n \ar[d] \\
\exp_kS^n \ar@{^{(}->}[r] & \exp_{k+1}S^n
}
\end{equation}
where $\Delta^{k+1}$ is the so-called {\em fat diagonal} in $SP^{k+1}(S^n)$ consisting of all those points where there is some duplication in coordinates. This has significantly more complicated structure than $(S^n)^{\times k}$, but if this were understood (together with the action of the Steenrod algebra on its cohomology), then that with the cohomology of the higher symmetric products $SP^k S^n$ should suffice to compute that of $Sub_k S^n$ for all $k$.

\bibliography{exp}

\end{document}